\documentclass[12pt, twoside]{article}
\usepackage{amsmath,amsthm}
\usepackage{amsfonts}
\usepackage{amssymb,latexsym}
\usepackage{enumerate}
\usepackage{cases}
\newtheorem{theorem}{Theorem}[section]
\newtheorem{remark}[theorem]{Remark}
\newtheorem{proposition}[theorem]{Proposition}
\newtheorem{lemma}[theorem]{Lemma}
\newtheorem{definition}[theorem]{Definition}

\newtheorem{corollary}[theorem]{Corollary}

\def\ind{1{\hskip -3 pt}\hbox{\textsc{I}}}
\def\n{\noindent}
\def\n{\noindent}

\def\Om{\Omega}
\def\E{\mathcal E}

\def\F{\mathcal F}

\def\n{\noindent}

\def\al{\alpha}

\def\R{\mathbb{R}}

\def\E{\mathcal E}
\def\N{\mathcal N}

\def\v{\hat{v}}
\def\u{\hat{u}}

\def\C{\mathbb{C}^n}

\def\O{\widehat{\Omega}}

\def\F{\mathcal F}
\DeclareMathOperator{\loc}{loc}
\numberwithin{equation}{section}
\usepackage{scrextend}
\changefontsizes[16pt]{13pt}
\begin{document}

\setlength{\baselineskip}{18truept}
\pagestyle{myheadings}
\markboth{  Van Phu Nguyen }{ Approximation of $m$-subharmonic functions  }
\title { Approximation of $m$-subharmonic functions in weighted energy classes with given boundary values}
\author{
	 Van Phu Nguyen
	\\ Faculty of Natural Sciences, Electric Power University,\\ Hanoi,Vietnam;
	\\E-mail: phunv@epu.edu.vn}

\date{}
\maketitle

\renewcommand{\thefootnote}{}



\renewcommand{\thefootnote}{\arabic{footnote}}
\setcounter{footnote}{0}

\begin{abstract}
This paper investigates the subextensions of m-subharmonic functions in weighted energy classes with given boundary values. These results are employed to approximate m-subharmonic functions in weighted energy classes by an increasing sequence of m-subharmonic functions defined on larger domains. In contrast to H. Amal's prior work with plurisubharmonic functions, our approach does not require the continuous boundary condition on the closure of the original domain.
\end{abstract}

\section{Introduction}

Subextension and approximation are pivotal concepts in pluripotential theory, with significant applications in complex geometry, especially in curvature problems such as prescribed Gauss curvature equations or Lagrangian mean curvature equations. Given domains $\Omega\subset \widehat{\Omega}$ be bounded domains in $\mathbb{C}^{n}$ and $u\in PSH(\Omega)$ (resp., $u\in SH_m(\Omega)$). A function $\hat{u}\in PSH(\widehat{\Omega})$ (resp., $\hat{u}\in SH_m(\widehat{\Omega})$) is said to be a subextension of $u$ if $\hat{u}\leq u$ on $\Omega.$ The concept of subextension in the class $\mathcal{F}(\Omega)$ has been studied by U. Cegrell and A. Zeriahi in 2003. In \cite{CZ03}, the authors demonstrated that if $\Omega\Subset \widehat{\Omega}$ are bounded hyperconvex domains in $\mathbb{C}^{n}$ and $u\in\mathcal{F}(\Omega)$ then there exists $\hat{u}\in\mathcal{F}(\widehat{\Omega})$ such that $\hat{u}\leq u$ on $\Omega$ and $\int_{\widehat{\Omega}}(dd^{c}\hat{u})^{n}\leq \int_{\Omega}(dd^{c}u)^{n}.$  The approximation in $\mathcal{F}(\Omega)$ has been studied by U. Cegrell and L. Hed in 2008. In \cite{CH08}, by applying the results of subextension, the authors proved that if $u\in\mathcal{F}(\Omega)$ and $u_j$ is the maximal subextention of $u$ to $\Omega_j$  then we have the following assertion: $u_j$ increases to $u$ a.e on $\Om$ is true for all $u\in\mathcal{F}(\Omega)$ if it is true for one function $u\in\mathcal{F}(\Omega).$ Approximation for weighted energy class $\mathcal{E}_{\chi}(\Omega)$ has been studied by S. Benekourch in \cite{Be11} and approximation with given boundary values in the class $\mathcal{F}(g|_{\Omega},\Om)$  has been studied by L. Hed  in \cite{H10}. After that, in 2014 , H. Amal \cite{A14} extended Hed's result from the class $\mathcal{F}(g|_{\Omega},\Om)$ to weighted energy class with given boundary $\mathcal{E}_{\chi}(g|_{\Omega},\Om).$\\
 Subextension and approximation for $m$-subharmonic function have been studied  recently by N. V. Phu and N. Q. Dieu in \cite{PDmmj}. In this paper, the authors proved a similar result to \cite{CH08} for the class of $m-$subharmonic functions (see Theorem 4.1) and extended also the main result in \cite{Be11} from plurisubharmonic functions to $m$-subharmonic functions (see Theorem 4.4). In \cite{Pjmaa} and \cite{AG23}, the authors generalized the main Theorem of  L. Hed (Theorem 1.1 in \cite{H10}) from the class of plurisubharmonic functions with given boundary values to the class of $m$-subharmonic functions with given boundary values (the two papers were conducted independently and yielded similar results). One of the novelties of the result in \cite{Pjmaa} and \cite{AG23} is that the author does not require the condition $g\in C(\overline{\Om})$.\\
 \n  Our goal is to extend these results from the case of plurisubharmonic functions to m-subharmonic functions in weighted energy classes with given boundary values. Firstly, we study subextension in weighted energy class with given boundary values of $m$-subharmonic functions (see Theorem \ref{th3.3}). Note that, compared to Theorem 4.1 in \cite{A14} in plurisubharmonic functions case, in Theorem \ref{th3.3}, we do not need the condition $f\in C(\overline{\Omega}).$ Secondly, using the result of subextension in weighted energy class with given boundary values of $m$-subharmonic functions, we study approximation in weighted energy class with given boundary values of $m$-subharmonic ones (see Theorem \ref{mainth}). The novelty of our result compared to H. Amal's main result in \cite{A14}  is that we do not need the condition $g\in C(\overline{\Omega}).$ To make our proof possible, we rely heavily on the method of constructing solution via upper envelopes of functions in suitable energy classes as was done in \cite{ACCH} and \cite{HHHP14}.
 This improvement is perhaps of some interest since in general a function in an energy class is only upper semicontinuous.\\
\n This paper is divided into four sections including the introduction and the remaining three parts. In Section 2 we refer readers to \cite{Pjmaa} for a summary of properities of  $m$-subharmonic functions which were studied intensively by many authors e.g. \cite{Bl1}, \cite{SA12}) \cite{Ch12}, \cite{T19} and \cite{Ch15}. We also recall some results on weighted $m$-enery classes in \cite{DT23} and some results that will frequently be used in the literature. In Section 3, we present some results related to subextension of the class $m$-subharmonic functions in weighted energy classes with given boundary values. In Section 4, we use the results in Section 3 to study the approximation of $m$-subharmonic functions in weighted energy classe with given boundary values. \\
{\bf Acknowledgments.}  
I wish to express our gratitude the Editor and  anonymous referees for their careful reading and constructive
comments that help to improve significantly my exposition. \\
The author would like to thank Professor Nguyen Quang Dieu for valuable comments during the preparation of this work.\\
The work is supported from Ministry of Education and Training, Vietnam under the Grant number B2025-CTT-10.\\
This work is revised in my visit in Vietnam Institute for Advanced Study in Mathematics (VIASM) in the Spring of 2025. I also thank VIASM for financial support and hospitality.	
\section{Preliminaries}
Some elements of the pluripotential theory that will be used throughout the paper can be found e.g. in \cite{BT1}, \cite{Ce98}, \cite{Ce04}, \cite{Kl}, while the elements of the theory of $m$-subharmonic functions and the complex $m$-Hessian operator can be found e.g. in \cite{Bl1}, \cite{Ch12}, \cite{Ch15}, \cite{SA12} and \cite{T19}. A summary of the properties required for this paper can be found in Preliminaries Section (from subsection 2.1 to subsection 2.6) in \cite{Pjmaa}.\\

\n Now we recall the definition of $\tilde{C}_m(E)-$capacity defined in \cite{HP17} as follows.
\begin{definition}
	Let $E$ be a Borel subset. We define
	$$\tilde{C}_m(E)= \tilde{C}_m(E,\Omega)= \sup\Bigl\{\int\limits_{E}H_{m-1}(u): u\in SH_m(\Omega), -1\leq u\leq 0\Bigl\}.$$
\end{definition}
\n According to the argument after Definition 2.7 in \cite{HP17}, we have $ \tilde{C}_m(E)\leq A(\Omega) {Cap}_m(E).$ It follows that if $u_j$ converges to $u$ in $Cap_m-$ capacity then $u_j$  also converges to $u$ in $\tilde{C}_m-$ capacity.\\
\n  We recall some results on weighted $m$-energy classes in \cite{DT23}. Let $\chi: \mathbb{R}^{-}\longrightarrow \mathbb{R}^{-},$ where $\mathbb{R}^{-}:= (-\infty, 0]$,
be an increasing function.
\begin{definition}{\rm We put
		\begin{align*}&\mathcal{E}_{m,\chi}(\Omega)=\{u\in SH_{m}(\Omega):\exists(u_{j})\in\mathcal{E}_{m}^{0}(\Omega), u_{j}\searrow u\,\,\text{and}\\
			& \hskip4cm \sup_{j}\int_{\Omega} -\chi(u_{j})H_{m}(u_{j})<+\infty\}.
	\end{align*}}
\end{definition}
\n Note that the weighted $m$-energy class generalizes the energy Cegrell class $\mathcal{F}_{m,p}$ and $\mathcal{F}_{m}.$
\begin{itemize}
	\item When $\chi\equiv -1,$ then $\mathcal{E}_{m,\chi}(\Omega)$ is the class $\mathcal{F}_{m}(\Omega).$\\
	\item When $\chi(t) =-(-t)^{p},$ then $\mathcal{E}_{m,\chi}(\Omega)$ is the class $\mathcal{E}_{m,p}(\Omega).$
\end{itemize}
According to Theorem 3.3 in \cite{DT23}, if $\chi\not\equiv 0$ then $\mathcal{E}_{m,\chi}(\Omega)\subset \mathcal{E}_{m}(\Omega)$ which means that the complex $m$-Hessian operator is well - defined on class $\mathcal{E}_{m,\chi}(\Omega)$ and if $\chi(-t)<0$ for all $t>0$ then we have $\mathcal{E}_{m,\chi}(\Omega)\subset \mathcal{N}_{m}(\Omega).$ 

\n We finish this Section with some results that will frequently be used in this paper.
\begin{theorem}[Theorem 3.8 in \cite{HP17}]\label{th3.8HP17}
	Let	$u_j,v_j,w\in\E_m(\Om)$ be such that $u_j,v_j\geq w,\forall j\geq 1.$ Assume that $|u_j-v_j|\to 0$ in $Cap_m$-capacity. Then we have $h[H_m(u_j)-H_m(v_j)]\to 0$ weakly as $j\to+\infty$ for all $h\in SH_m\cap L^{\infty}_{\loc}(\Om).$
\end{theorem}
\n Lemma 2.9 in \cite{T19} stated that if $[u_j]$ is a monotone sequence of $m$-subharmonic functions converging to a $m$-subharmonic function $u$ then we have $u_j\to u$ in $Cap_m$ as $j\to\infty.$ Thus, we have the following Corollary.
\begin{corollary}\label{co3.8HP17}
	Let $[u_j]\subset\E_m(\Om)$ be a monotone sequence of functions converging to a function $u \in\E_m(\Om)$ then the sequence of measures $H_m(u_j)\to H_m(u)$ weakly as $j\to\infty.$
\end{corollary}
\begin{proposition}[Proposition 2.9 in \cite{Pjmaa}]\label{pro5.2HP17}
	Let $u,v,u_k\in\E_m(\Om), k=1,\cdots,m-1$ with $u\geq v $ on $\Om$ and  $T=dd^cu_1\wedge\cdots\wedge dd^cu_{m-1}\wedge\beta^{n-m}.$ Then we have
	$$\ind_{\{u=-\infty\}}dd^cu\wedge T\leq \ind_{\{v=-\infty\}}dd^cv\wedge T.$$
	In particular, if $u,v\in\E_m(\Om)$ are such that $u\geq v$ then for every $m-$polar set $A\subset\Om$ we have
	$$ \int_{A}H_m(u)\leq\int_{A}H_m(v).$$	
\end{proposition}

\begin{lemma}[Lemma 2.14 in \cite{Gasmi}]\label{lm2.7Gasmi}
	Let $f\in\mathcal{E}_m(\Omega)$ and $u,v\in\mathcal{N}_m(f)$ are such that $u\leq v$ and $\int_{\Omega}H_m(u)<+\infty.$ Then for all $\varphi\in SH_m^{-}(\Om)$ the following inequality holds
	$$\int_{\Om}H_m(v)\leq\int_{\Om}H_m(u).$$
\end{lemma}
\begin{lemma} [Theorem 3.1 in \cite{Pjmaa}]\label{l3.1}
	If $u\in\mathcal{N}_m(f)$ satisfies $\int_{\Om}H_m(u)<+\infty$ then $u\in\F_m(f).$
\end{lemma}

\section{Subextension of $m$-subharmonic functions in weighted energy classes with given boundary values}
\n Firstly, we prove the following key Lemma which is inspired from Lemma 4.13 in \cite{ACCH} (see also \cite{HHHP14}).

\begin{lemma}\label{lm3.2}
	Let $\Om$ be a bounded $m$-hyperconvex domain in $\C$ and let $\alpha$ be a finite Radon measure on $\Om$ which puts no mass on $m$-polar sets of $\Om$ such that supp$\alpha\Subset\Om.$ Assume that $\chi: \mathbb{R}^-\to\mathbb{R}^-$ is a nondecreasing continuous function such that $\chi(t)<0$ for all $t<0,$ $\chi(-\infty)=-1$ and $f\in\E_m(\Om)\cap MSH_m(\Om).$ Let $v\in \mathcal F_{m}(f,\Omega)$  be such that $\text{\rm supp}H_{m}(v)\Subset\Omega,$  $H_{m}(v)$ is carried by a $m-$polar set and $\int_{\Om}H_m(v)<+\infty.$ Set 
	$$\mathcal U(\alpha, v)=\left\{\varphi\in SH^{-}_{m}(\Omega): \alpha\leq -\chi (\varphi)H_{m}(\varphi) \text{ and } \varphi\leq v \right\}.$$
	We define
	$$u=(\sup \left\{\varphi: \varphi\in \mathcal U(\alpha, v)\right\})^*.$$
	Then we have $u\in\F_m(f,\Om)$ and
	$-\chi(u)H_m(u)=\alpha+H_m(v).$
\end{lemma}
\begin{proof}
	We split the proof into two steps as follows.
	
	\n 	
	{\em Step 1.}
	We prove that $u\in\mathcal{F}_m(f,\Omega)$ and
	\begin{equation}\label{eq5}
		-\chi (u)H_{m}(u)\geq\alpha -\chi(v)H_{m}(v).
	\end{equation}
	Indeed, first by Theorem 4.1 in \cite{PD2023} there exists a function $\phi\in\mathcal{E}_{m,\chi}(\Om)\subset \mathcal N_{m}(\Omega)$ such that 
	$$ -\chi(\phi)H_{m}(\phi)=\alpha.$$
	\n
	Note that supp$H_{m}(\phi) \Subset \Omega$ so $\int_\Omega H_{m}(\phi)<+\infty$. Thus, by Theorem 4.9 in \cite{T19} we infer that $\phi\in\mathcal{F}_m(\Omega).$
	Moreover, it is easy to see that $(\phi+v)$ belongs to $\mathcal U(\alpha,v).$ So we have  $$\phi+v\leq u\leq v.$$ It follows from $v\in\F_m(f,\Om)$ that $u\in\mathcal{F}_m(f,\Omega).$\\
	Next, using Choquet Lemma we may find a sequence of functions $[\xi_j]\subset \mathcal U(\alpha,v)$ such that
	$$u=\Bigl(\sup_{j\in \mathbb N^*} \xi_j\Bigl)^*.$$
	Note that, by Lemma 3.1 in \cite{HQ21} we have $\max (\varphi,\psi)\in\mathcal U(\alpha, v)$, $\forall\ \varphi,\psi\in \mathcal U(\alpha, v)$. Thus, 
	$$\tilde u_j:=\max\{\xi_1,\ldots,\xi_j\} \in \mathcal U(\alpha,v)$$ and $\tilde u_j\nearrow u$ in $Cap_m$-capacity. By Corollary 3.3 in \cite{PD2023} we also have $$-\chi (\tilde u_j)H_{m}(\tilde{u_{j}})\to -\chi (u)H_{m}(u)$$ weakly. Therefore it follows that $$-\chi (u)H_{m}(u)\geq\alpha.$$ Thus, we get  $u\in\mathcal U(\alpha,v)$. Now, since $\alpha$ vanishes on $m$-polar sets, we have
	\begin{align}\label{e3.11}
		\begin{split}
			-\chi(u)H_{m}(u) &= -\chi(u) \ind_{\{u>-\infty\}} H_{m}(u) -\chi(u) \ind_{\{u=-\infty\}} H_{m}(u)
			\\
			&\geq  \ind_{\{u>-\infty\}} \al -\chi(u) \ind_{\{u=-\infty\}} H_{m}(u)\\
			&\geq \al -\chi(u) \ind_{\{u=-\infty\}} H_{m}(u).
		\end{split}
	\end{align}
	On the other hand, since $H_{m}(v)$ is carried by an $m$-polar set then we have
	$$H_{m}(v)=\ind_{\{v>-\infty\}}H_{m}(v)+ \ind_{\{v=-\infty\}}H_{m}(v)=\ind_{\{v=-\infty\}}H_{m}(v).$$
	It follows from $u\leq v$ and Proposition \ref{pro5.2HP17} that
	\begin{equation}\label{e3.12} 1_{\{u=-\infty\}} H_{m}(u)\geq 1_{\{v=-\infty\}} H_{m}(v) =H_{m}(v).
	\end{equation}
	Combining inequalities \ref{e3.11} and \ref{e3.12} we get
	$$-\chi(u)H_m(u)\geq \alpha -\chi (u) H_{m}(v)\geq \alpha -\chi (v) H_{m}(v).$$
	
	{\em Step 2.} We prove that
	$$-\chi(u)H_m(u) = \alpha -\chi (v) H_{m}(v).$$
	Choose a $m$-hyperconvex domain $G\Subset\Om$ such that $\text{\rm supp}\alpha \cup\text{\rm supp}H_{m}(v)\Subset G\Subset \Om$. Next, since $v\in \mathcal{F}_{m}(f,\Om)$ and Lemma 2.13 in \cite{PD2024}, we can choose a sequence
	$\mathcal E_{m}^{0}(f,\Omega) \ni \{v_j\}\searrow v $ such that $ \text{\rm supp}H_{m}(v_{j})\subset\overline{G}$ and
 $H_m(v_j)$ vanishes on all $m$-polar sets.
	Thus, by Theorem  3.1 in \cite{PD2024}, there exists a unique function $w_j\in\mathcal E_{m,\chi}(f,\Omega)$ such that
	\begin{equation} \label{eq7}
		-\chi(w_j) H_{m}(w_{j}) =\alpha - \chi( v_j )H_{m}(v_{j}).\end{equation}
	On the other hand, we have
	$$ -\chi(\phi+v_j)H_m(\phi)\geq -\chi(\phi)H_m(\phi)=\alpha$$
	and 
	$$ -\chi(\phi+v_j)H_m(v_j)\geq -\chi(v_j)H_m(v_j)$$
	Thus, we have 
	\begin{align*}-\chi (\phi+v_j) H_m (\phi+v_j)&\geq\alpha -\chi(v_j)H_m(v_j)\\
		&= -\chi(w_j) H_{m}(w_j)\\
		&\geq -\chi (v_j) H_{m}(v_j).
	\end{align*} 
	By Theorem 3.8 in \cite{PDjmaa} we get $$\phi+v_j\leq w_j\leq v_j.$$ It follows from $\phi\in\F_m(\Om)$ and $v_j\in\mathcal{E}_m^0(f,\Om)$ that $w_j\in\mathcal{F}_m(f,\Omega).$ We put 
	$$\mathcal V(\alpha, v_j)=\left\{\varphi\in SH^{-}_{m}(\Omega): \alpha\leq -\chi (\varphi)H_{m}(\varphi) \text{ and } \varphi\leq v_j \right\}.$$ Obviously, we have $w_j\in\mathcal V(\alpha, v_j).$ 
	
	Now we define
	$$u_j:=\bigl(\sup \left\{\varphi: \varphi\in \mathcal {V}(\alpha,v_j)\right\}\bigl)^{*}.$$
	Then $u_j$ is decreasing and $u_j\geq w_j$ for all $j\geq 1$. Note that, by $\mathcal{F}_m(f,\Omega)\ni w_j\leq u_j,$ we also have $u_j\in\mathcal{F}_m(f,\Omega).$ We will prove that $\{u_j\}$  decreases pointwise to $u$, as $j\to+\infty$. Indeed, assume that $u_j\searrow \tilde{u}$ as $j\to +\infty.$\\
	\n (i) We have $u\leq v\leq v_j.$ Thus, $u\in \mathcal{V}(\alpha,v_j).$ It follows that $u\leq u_j.$ Let $j\to+\infty$ we get that $u\leq \tilde{u}.$\\
	\n (ii) On the other hand, we have $u_j\leq v_j.$ Let $j\to+\infty$ we obtain that $\tilde{u}\leq v.$ Note that, since $u_j\searrow \tilde{u}$ we get that $u_j$ converges to $\tilde{u}$ in $Cap_m-$ capacity. By Corollary 3.3 in \cite{PD2023} we get that $-\chi(u_j)H_m(u_j) \to -\chi(\tilde{u})H_m(\tilde{u}) \,\,\text{weakly}.$
	Therefore $-\chi(\tilde{u})H_m(\tilde{u})\geq \alpha.$ So we have $\tilde{u}\in\mathcal{U}(\alpha,v)$ and we obtain $\tilde{u}\leq u.$\\ 
	From (i) and (ii) we get $\tilde{u}=u.$ That means	$u_j$  decreases pointwise to $u$, as $j\to+\infty$.\\
	Next,  we will prove that $u_j-w_j\to 0$ in $\tilde{C}_{m}$-capacity. Indeed, fix a strictly $m$-subharmonic function $h_0\in \mathcal E_m^{0}(\Om)\cap \mathcal C^\infty(\Om)$, such a function $h_0$ exists by \cite{ACL18}. For  $\varepsilon>0$ and $j_0\ge 1$ we have
	\begin{align*}
		& \limsup_{j\to\infty} \tilde{C}_{m} (\{u_j-w_j>\varepsilon\})
		\\&=\limsup_{j\to\infty}\Bigl( \sup\Bigl\{\int_{\{u_j-w_j>\varepsilon\}} dd^c h_0\wedge (dd^c h)^{m-1}\wedge\beta^{n-m}:
		\\&\hspace{9cm} h\in SH_m(\Omega), -1\leq h\leq 0\Bigl\}\Bigl)
		\\&\leq\limsup_{j\to\infty}\Bigl( \sup\Bigl\{\frac{1}{ \varepsilon^m }\int_{\{u_j-w_j>\varepsilon\}} (u_j-w_j)^m  dd^c h_0\wedge (dd^c h)^{m-1}\wedge\beta^{n-m}:
		\\&\hspace{9cm} h\in SH_m(\Omega), -1\leq h\leq 0\Bigl\}\Bigl)
		\\&\leq\limsup_{j\to\infty}\Bigl( \sup\Bigl\{\frac{1}{\varepsilon^m }\int_{\Omega} (u_j-w_j)^m  dd^c h_0\wedge (dd^c h)^{m-1}\wedge\beta^{n-m}:
		\\&\hspace{9cm} h\in SH_m(\Omega), -1\leq h\leq 0\Bigl\}\Bigl)
\\&\leq\limsup_{j\to\infty} \frac{m!}{\varepsilon^m} \int_{\Omega} -h_0 \Bigl[H_m(w_j)-H_m(u_j)\Bigl]
\\&\leq\limsup_{j\to\infty} \frac{m!}{ \varepsilon^m } \int_{\Omega} -h_0 H_m(w_j) +\frac {m!} { \varepsilon^m } \int_{\Omega} h_0 H_{m}(u)
\\&\leq\limsup_{j\to\infty} \frac {m!} { \varepsilon^m } \int_{\Omega} -h_0 \frac {-\chi (w_j) H_m(w_j)} {-\chi (u_j)} +\frac {m!} { \varepsilon^m } \int_{\Omega} h_0 H_m(u)
\\&\leq\limsup_{j\to\infty} \frac {m!} { \varepsilon^m } \int_{\Omega} -h_0 \frac {\alpha -\chi(v_j) H_m(v_j)} {-\chi (u_j)} +\frac {m!} { \varepsilon^m } \int_{\Omega} h_0 H_m(u)
\\&\leq\limsup_{j\to\infty} \frac {m!} { \varepsilon^m } \int_{\Omega} -h_0 \frac {\alpha + H_m(v_j)} {-\chi (u_j)} +\frac {m!} { \varepsilon^m } \int_{\Omega} h_0 H_m(u)
\\&\leq\limsup_{j\to\infty} \frac {m!} { \varepsilon^m } \int_{\Omega} -h_0 \frac {\alpha + H_m(v_j)} {-\chi (u_{j_0})} +\frac {m!} { \varepsilon^m } \int_{\Omega} h_0 H_m(u),
\end{align*}
	where the fourth estimate follows from Lemma 5.4 in \cite{T19} and the sixth estimate follows from  supp$H_m(w_j)\subset \overline{G}\Subset\Omega$ and $\chi(t)<0$ for all $t<0$ which imply that $\chi(u_j)<0.$
	Since $v_j\searrow v,$ we infer that $v_j\to v$ in $Cap_m$-capacity. Thus, it follows from Theorem 3.8 in \cite{HP17} that $H_m(v_j)\to H_m(v)$ weakly as $j\to+\infty.$ Moreover, $\chi(u_{j_{0}})$ is an upper semicontinuous function, so by Lemma 1.9 in \cite{De93} we have 
	$$\lim\limits_{j\to\infty}\frac { H_m(v_j)} {-\chi (u_{j_0})}\leq \frac { H_m(v)} {-\chi (u_{j_0})}.$$
	\n Note that $\text{\rm supp}\alpha \cup\text{\rm supp}H_m(v_j)\subset\overline{G}$, we obtain
	\begin{align*}
		&\limsup_{j\to\infty}  \int_{\Omega} -h_0 \frac {\alpha + H_m(v_j)} {-\chi (u_{j_0})} \\
		&\leq \limsup_{j\to\infty}  \int_{\overline{G}} -h_0 \frac {\alpha + H_m(v_j)} {-\chi (u_{j_0})} \\
		&\leq  \int_{\overline{G}} -h_0 \frac {\alpha + H_{m}(v)} {-\chi (u_{j_0})} \\
		&\leq  \int_{\Omega} -h_0 \frac {\alpha + H_{m}(v)} {-\chi (u_{j_0})}. 
	\end{align*}
	Thus, we have
	\begin{align*}
		& \limsup_{j\to\infty} \tilde{C}_{m} (\{u_j-w_j>\varepsilon\})\\
		&\leq \frac {m!} { \varepsilon^m } \int_{\Omega} -h_0 \frac {\alpha + H_m(v)} {-\chi (u_{j_0})} +\frac {m!} { \varepsilon^m } \int_{\Omega} h_0 H_m(u).
	\end{align*}	
	
	\n Letting $j_0\to\infty$ by the Lebegue monotone convergence Theorem and Step 1 we get
	\begin{align*}
		\limsup_{j\to\infty} \tilde{C}_{m} (\{u_j-w_j>\varepsilon\})&\leq\limsup_{j_0\to\infty} \frac {m!} { \varepsilon^m } \int_{\Omega} -h_0 \frac {\alpha + H_m(v)} {-\chi (u_{j_0})} +\frac {m!} { \varepsilon^m } \int_{\Omega} h_0 H_m(u)
		\\&= \frac {m!} { \varepsilon^m } \int_{\Omega} -h_0 [\frac {\alpha + H_m(v)} {-\chi (u)} - H_m(u)]\leq 0.
	\end{align*}
\n Thus we have proved that $u_j-w_j \to 0$ in $\tilde{C}_{m}$-capacity.
	Since $u_j \to u$ in $Cap_m$-capacity, we infer that $w_j\to u$ in $\tilde{C}_{m}$-capacity. Hence by Proposition 3.1 in \cite{PD2023} we get
	$$\liminf\limits_{j\to\infty} -\chi(w_j) H_m(w_j)\geq -\chi(u) H_m(u).$$
	This is equivalent to 
	$$\liminf\limits_{j\to\infty}[ \alpha -\chi(v_j) H_m(v_j)]\geq -\chi(u) H_m(u).$$
	Note that it follows from $v_j\searrow v$ and Lemma 2.9 in \cite{T19} that $v_j\to v$ in $Cap_m$-capacity. By Corollary 3.3 in \cite{PD2023} we have 
	$$\lim\limits_{j\to\infty} -\chi(v_j) H_m(v_j)=-\chi(v) H_m(v).$$
	This implies that
	 $$\alpha -\chi(v)H_m(v)\geq -\chi(u) H_m(u).$$
	Combining this with inequality \eqref{eq5}, we finally reach
	$$-\chi(u) H_m(u)=\alpha -\chi(v) H_m(v).$$
\end{proof}
\begin{remark}
	
	In Step 1 in the proof of the Lemma $\ref{lm3.2}$, if we omit the hypothesis supp$\alpha\Subset\Om$ then we only have $u\in\mathcal{N}_{m}(f,\Om).$ 
\end{remark}

\n Now we are ready for  the main result of this Section about subextension in the class $\mathcal{E}_{m,\chi}(f)$. Note that, compared to Theorem 4.1 in \cite{A14} in plurisubharmonic functions case, in the theorem below, we do not have the condition $f\in C(\overline{\Omega}).$
\begin{theorem}\label{th3.3}
	Let $\Om\Subset\O$ are bounded $m$-hyperconvex domains in $\C$ and $f\in\E_m\cap MSH_m(\Om), g\in\E_m(\O)\cap MSH_m(\O)$ with $f\geq g$ on $\Om.$ Assume that $\chi: \mathbb{R}^-\to\mathbb{R}^-$ be a nondecreasing continuous function such that $\chi(t)<0$ for all $t<0.$ Then for every $u\in\mathcal{E}_{m,\chi}(f,\Om)$ such that 
	\begin{equation}\label{eq3.8a}\int_{\Om}[-\chi(u)+1]H_m(u)<+\infty,
	\end{equation}
	there exists $\u\in\E_{m,\chi}(g,\O)$ such that $\u\leq u$ on $\Om$ and 
	$$-\chi(\u)H_m(\u)=\ind_{\Om}[-\chi(u)H_m(u)]\,\,\text{on}\,\, \O.$$
	
\end{theorem}
\begin{proof} 
	We outline some main points of our proof. First, we find a lower bound for $u$ in the class $\mathcal{E}_{m}(\Om).$ Then the proof is divided into two cases $\chi(-\infty)>-\infty$ and $\chi(-\infty)=-\infty.$ In the first case, we will construct in  Step 3 the function $\hat{u}$ as an upper envelope of a class of functions in the energy class $\mathcal{E}_m(\hat{\Omega}).$ This class of functions will be defined in  Step 1 and Step 2.
	  In the remaining Step 4 and Step 5, we will prove that $\hat{u}$ satisfies the condition of Theorem. For the Case 2, the existence of $\hat{u}$ is directly deduced from a recent work of \cite{PD2024}. Now, we give a detail proof.\\
	Since $u\in\mathcal{E}_{m,\chi}(f,\Om),$ there exists a function $\varphi\in\mathcal{E}_{m,\chi}(\Om)$ such that $f\geq u\geq f+\varphi.$ According to the assumption $\chi(t)<0$ for all $t<0$ and Theorem 3.3 in \cite{DT23}, we obtain that $\varphi\in \N_m(\Om).$\\
	We split the proof into two cases.\\
	\n {\bf Case 1.} When $\chi(-\infty)>-\infty.$ After multiplying with a constant, we can assume that $\chi(-\infty)=-1.$\\
\n	We split the proof of Case 1 into five steps as below.\\
\n {\em Step 1.} We will show that there exist	$v\in\F_m(f,\Om)$ such that $v\geq u$ and 
$$H_m(v)=\ind_{\{u=-\infty\}}H_m(u).$$
	Indeed, it follows from the assumption \eqref{eq3.8a} that we have
	\begin{equation}\label{e3.1}\int_{\Om\cap\{u=-\infty\}}H_m(u)=\int_{\Om\cap\{u=-\infty\}}-\chi(u)H_m(u)\leq\int_{\Om}-\chi(u)H_m(u)<+\infty.\end{equation}
 We put $\tau=\max(\varphi,u),$ then we have $\tau\in\N_m(\Om)$ and $\tau\geq u$ on $\Om.$ Obviously, we have $f\geq u\geq f+\tau.$ Thus, by replacing $\varphi$ by $\tau,$ we can assume that $\varphi\geq u.$ Hence, we have $0\leq \varphi-u\leq-f.$ Note that $f\in MSH_m(\Om),$ by Lemma 5.1 in \cite{Gasmi}, we infer that
$$\ind_{\{u=-\infty\}}H_m(u)= \ind_{\{\varphi=-\infty\}}H_m(\varphi)\leq H_m(\varphi).$$
According to the main Theorem in \cite{Gasmi},  there exists a function $\kappa\in\N_m(f,\Om)$ such that $$\ind_{\{u=-\infty\}}H_m(u)=H_m(\kappa)$$ and
\begin{equation}\label{e3.2} f\geq \kappa\geq f+\varphi.
\end{equation}
Now, by the inequality $\eqref{e3.1}$, we also have $$\int_{\Om}H_m(\kappa)= \int_{\Om}\ind_{\{u=-\infty\}}H_m(u)<+\infty.$$ It follows from Lemma $\ref{l3.1}$ that $\kappa\in\F_m(f,\Om).$\\
On the other hand, by Theorem 4.8 in \cite{T19}, we obtain $\kappa\in\F_m(\tilde{\kappa}),$ where $\tilde{\kappa}$ is defined as in the Subsection 2.6 in \cite{Pjmaa}. In detail, there exists a function $\kappa_1\in \F_m(\Om)$ such that $$\kappa\geq \tilde{\kappa}+\kappa_1$$ and 
\begin{equation}\label{eq3.1} \int_{\Om}H_m(\kappa_1)\leq\int_{\Om}H_m(\kappa).
\end{equation} 
Note that, it follows from $\varphi\in\mathcal{N}_m(\Omega)$ that $\tilde{\varphi}=0.$ Thus, from inequality \ref{e3.2}, we infer that $\tilde{\kappa}=\tilde{f}.$ Therefore, we have $$\kappa\geq \tilde{\kappa}+\kappa_1= \tilde{f}+\kappa_1\geq f+\kappa_1.$$
We put 
$$\rho=(\sup\{\zeta\in SH_m^-(\Om):\zeta + f\leq \kappa\,\,\text{on}\,\,\Om\})^*.$$
Then we have  $\F_m(\Om)\ni \kappa_1\leq  \rho.$ By integration by parts, we have
\begin{equation}\label{eq3.2}
	\int_{\Om}H_m(\rho)\leq\int_{\Om}H_m(\kappa_1).
\end{equation}
Combining inequalities $\eqref{eq3.1}$ and $\eqref{eq3.2},$ we have
\begin{equation}\label{eq3.3}
	\int_{\Om}H_m(\rho)\leq \int_{\Om}H_m(\kappa).
\end{equation}
Moreover, from $\rho +f\leq \kappa$ on $\Om$ and Proposition $\ref{pro5.2HP17},$ we get 
$$\ind_{\{\kappa=-\infty\}}H_m(\kappa)\leq \ind_{\{\rho+f=-\infty\}}H_m(\rho+f).$$
Note that $f\in MSH_m(\Om)$ and $|(f+\rho)-\rho|=-f.$ From Lemma 5.1 in \cite{Gasmi}, we have
$$\ind_{\{\rho+f=-\infty\}}H_m(\rho+f)=\ind_{\{\rho=-\infty\}}H_m(\rho).$$
Thus, we obtain
$$\ind_{\{\kappa=-\infty\}}H_m(\kappa)\leq \ind_{\{\rho=-\infty\}}H_m(\rho).$$
Moreover, since $H_m(\kappa)$ is carried by a $m$-polar set, so we have
	$$	H_m(\kappa)=\ind_{\{\kappa=-\infty\}}H_m(\kappa)
		\leq \ind_{\{\rho=-\infty\}}H_m(\rho)\\
		\leq H_m(\rho).$$
Combining with inequalities $\eqref{eq3.3}$  we deduce that 
$$ H_m(\kappa)=\ind_{\{\rho=-\infty\}}H_m(\rho).$$
This is equivalent to
\begin{equation}\label{eq3.5}
	\ind_{\{u=-\infty\}}H_m(u)=\ind_{\{\rho=-\infty\}}H_m(\rho).
\end{equation}
Now, we put $v=\max(\kappa,u).$ Obviously, we have $v\in\mathcal{N}_m(f,\Om)$ and $v\geq u.$ We claim that $$H_m(v)=\ind_{\{u=-\infty\}}H_m(u).$$
Indeed, it follows from inequality \ref{e3.2} that $\kappa\geq f+v.$ Therefore, by the construction of $\rho,$ we obtain $u\leq v\leq \rho.$ Thus, by Proposition $\ref{pro5.2HP17}$ we have
$$\ind_{\{\rho=-\infty\}}H_m(\rho)\leq \ind_{\{v=-\infty\}}H_m(v)\leq\ind_{\{u=-\infty\}}H_m(u).$$
Equation $\eqref{eq3.5}$ implies that
$$\ind_{\{\rho=-\infty\}}H_m(\rho)= \ind_{\{v=-\infty\}}H_m(v)=\ind_{\{u=-\infty\}}H_m(u).$$
Therefore, we have
\begin{equation}\label{eq3.6}
	H_m(\kappa)=\ind_{\{u=-\infty\}}H_m(u)=\ind_{\{v=-\infty\}}H_m(v)\leq H_m(v).
\end{equation}
Moreover, since $v,\kappa\in \N_m(f,\Om), v\geq \kappa$ and $\int_{\Om}H_m(\kappa)<+\infty,$ according to Lemma $\ref{lm2.7Gasmi}$ we have
\begin{equation}\label{eq3.7}
	\int_{\Om}H_m(v)\leq 	\int_{\Om}H_m(\kappa).
\end{equation}
Combining inequation $\eqref{eq3.6}$ and inequation $\eqref{eq3.7}$ we have
$$H_m(v)=H_m(\kappa)=\ind_{\{u=-\infty\}}H_m(u).$$ 
Moreover, we have
$\int_{\Om}H_m(v)=\int_{\Om}\ind_{\{u=-\infty\}}H_m(u)<+\infty$ and $v\in\N_m(f,\Om).$ By Lemma $\ref{l3.1}$ we infer that $v\in\F_m(f,\Om).$ The proof of Step 1 is complete.\\

\n {\em Step 2.} Put 
	$$\hat{v}=(\sup\{\xi\in\mathcal{E}_m(g,\O):\xi\leq v\,\,\text{on}\,\,\Om\})^*.$$
	\n We will prove that $\hat{v}\in\F_m(g,\O),\v\leq v$ on $\Om$ and $H_m(\v)=\ind_{\Om}H_m(v).$ Indeed, by Theorem 4.4 in \cite{Pjmaa} there exist $\phi\in\F_m(g,\O)$ such that $\phi\leq v$ on $\Omega$ and $H_m(\phi)\leq\ind_{\Om}H_m(v)$ on $\O.$ This implies that
	\begin{equation}\label{eq3.8}
		\int_{\O}H_m(\phi)\leq\int_{\Om}H_m(v)<+\infty.
	\end{equation}
	On the other hand, we have $\phi$ belongs to the class of functions which are in the definition of $\v.$ Therefore, we infer that $\phi\leq\v$ on $\O.$ It follows from Lemma $\ref{lm2.7Gasmi}$ that
	\begin{equation}\label{eq3.9}
		\int_{\O}H_m(\v)\leq\int_{\O}H_m(\phi).
	\end{equation}
	Combining inequality $\eqref{eq3.8}$ and inequality $\eqref{eq3.9}$ we obtain
	\begin{equation}\label{eq3.10}
		\int_{\O}H_m(\v)\leq\int_{\Om}H_m(v)=\int_{\O}\ind_{\Om}H_m(v).
	\end{equation}
	Moreover, since $\v\leq v$ on $\Om$ and Proposition $\ref{pro5.2HP17}$ we infer that
	$$\ind_{\{v=-\infty\}}H_m(v)\leq \ind_{\{\v=-\infty\}}H_m(\v)\,\,\text{on}\,\,\Om.$$
	Note that $H_m(v)$ is carried by an $m$-polar set. This means that $H_m(v)= \ind_{\{v=-\infty\}}H_m(v).$ Therefore, we have
	\begin{equation}\label{eq3.11}
		\ind_{\Om}H_m(v)\leq \ind_{\{\v=-\infty\}}H_m(\v)\leq H_m(\v).
	\end{equation}
	Combining inequality $\eqref{eq3.10}$ and inequality $\eqref{eq3.11}$, we obtain that
	$$H_m(\v)=\ind_{\Om}H_m(v).$$
It follows from $\v\geq\phi$ and $\phi\in\mathcal{F}_m(g,\O)$ that $\v\in\mathcal{F}_m(g,\O).$	The proof of Step 2 is complete.\\

\n {\em Step 3.} We will show that there exists $\hat{u}\in\mathcal{F}_{m}(g,\O)$ such that  
$$-\chi(\u)H_m(\u)=\ind_{\Om}[-\chi(u)H_m(u)]\,\,\text{on}\,\, \O.$$
Indeed, we put 
	$$\u=(\sup \left\{\xi: \xi\in \mathcal U(\tilde{\alpha}, \v)\right\})^*,$$
	where
	$$\mathcal U(\tilde{\alpha}, \v)=\left\{\xi\in \mathcal E_{m}(\O): \tilde{\alpha}\leq -\chi (\xi)H_{m}(\xi) \text{ and } \xi\leq \v \right\}$$ with $\tilde{\alpha}=\ind_{\{u>-\infty\}\cap\Om}[-\chi(u)]H_m(u)$ is a nonnegative measure that vanishes on all $m$-polar sets and supp$H_m(\tilde{\alpha})\Subset\O.$ Note that,  we also have $H_m(\v)$ is carried by a $m$-polar set and supp$H_m(\v)\Subset\O.$
	According to Lemma $\ref{lm3.2}$ we get $\u\in\F_m(g,\O)$ and
	\begin{align}\label{eq3.12b}
		\begin{split}
			-\chi(\u)H_m(\u)&=\tilde{\alpha}+H_m(\v)\\
			&=\ind_{\{u>-\infty\}\cap\Om}[-\chi(u)]H_m(u)+ \ind_{\Om}H_m(v)\\
			&= \ind_{\{u>-\infty\}\cap\Om}[-\chi(u)]H_m(u)+ \ind_{\Om}\ind_{\{u=-\infty\}}H_m(u)\\
			&=\ind_{\{u>-\infty\}\cap\Om}[-\chi(u)]H_m(u)+\ind_{\{u=-\infty\}\cap\Om}[-\chi(u)]H_m(u)\\
			&=\ind_{\Om}[-\chi(u)]H_m(u).
	\end{split}
	\end{align} 
\n The proof of Step 3 is complete.\\

\n {\em Step 4.} We will prove that $\hat{u}\leq u$ on $\Omega.$

Indeed, we put 
$$w=(\sup \left\{\xi: \xi\in \mathcal U(\alpha, v)\right\})^*,$$
where
$$\mathcal U(\alpha, v)=\left\{\xi\in \mathcal E_{m}(\Omega): \alpha\leq -\chi (\xi)H_{m}(\xi) \text{ and } \xi\leq v \right\}$$ with $\alpha=\ind_{\{u>-\infty\}}[-\chi(u)]H_m(u)$ is a nonnegative measure that vanishes on all $m$-polar sets and $H_m(v)$ is carried by a $m$-polar sets.\\
It is easy to check that $u$ belongs to the class of functions which are in the definition of $w.$ This implies that $u\leq w.$\\
\n Repeating argument as in Step 1 in the proof of Lemma $\ref{lm3.2}$ with note that without the assumption supp$\alpha\Subset\Om$ we only get $w\in\N_m(f,\Om)$ and
\begin{align*} -\chi(w)H_m(w)&\geq \alpha+ H_m(v)\\
	&\geq \ind_{\{u>-\infty\}}[-\chi(u)]H_m(u)+\ind_{\{u=-\infty\}}H_m(u)\\
	&=\ind_{\{u>-\infty\}}[-\chi(u)]H_m(u)+\ind_{\{u=-\infty\}}[-\chi(u)]H_m(u)\\
	&\geq -\chi(u)H_m(u)\\
	&\geq -\chi(w)H_m(u).
\end{align*}
This implies that $H_m(w)\geq H_m(u)$ on $\Om.$ Moreover, it follows from $u\leq w,$  $\int_{\Om}H_m(u)<+\infty$ and Lemma $\ref{lm2.7Gasmi}$ that $\int_{\Om}H_m(w)\leq\int_{\Om}H_m(u).$ Thus, we obtain that $H_m(w)=H_m(u).$ By Theorem 2.10 in \cite{Gasmi2} (see Theorem 3.6 in \cite{ACCH} for the case of plurisubharmonic functions), we obtain that $w=u.$ By the constructions of $u$ and $\u$ we get $\u\leq u$ on $\Om.$ The proof of Step 4 is complete.\\

\n {\em Step 5.}	
	 It remains to prove that $\u\in\mathcal{E}_{m,\chi}(g,\O).$
	We claim that $\u\in MSH_m(\O\smallsetminus\overline{\Om}).$ Indeed, let $\psi\in SH_m(\O\smallsetminus\overline{\Om})$ and $K\Subset \O\smallsetminus\overline{\Om}$ such that $\psi\leq\u$ on $\partial K.$ We put
	$$	\gamma=
	\begin{cases}
		\max(\psi,\u) & \,\,\text{on}\,\, K\\
		\u&\,\,\text{on}\,\,\O\smallsetminus K.	
	\end{cases}
	$$
	Then we have $\gamma\in SH_m(\O).$ We will prove that $\gamma$ belongs to the class of functions which are in the definition of $\u.$\\
	Firstly, we prove that 
	\begin{equation}\label{eq3.12a}
		\gamma\leq\v.
	\end{equation}
	Indeed, it follows from $\gamma\leq g$ on $\partial K$ and $g\in MSH_m(\O)$ that $\gamma\leq g$ on $\O.$ Moreover, we also have $\gamma\geq\u\in\F_m(g,\O).$ Thus, we get $\gamma\in\mathcal{F}_m(g,\O).$ Now, on $\Om$ we have
	$$\gamma=\u\leq\v\leq v.$$
	Thus, by the construction of $\v,$ we have $\gamma\leq\v$ on $\O$ and we have inequality $\eqref{eq3.12a}.$\\
	Secondly, we will prove that 
	\begin{equation}\label{eq3.15}
		-\chi(\gamma)H_m(\gamma)\geq  \tilde{\alpha}=\ind_{\{u>-\infty\}\cap\Om}[-\chi(u)]H_m(u).
	\end{equation}
	Indeed, since $$-\chi(\u)H_m(\u)\geq \ind_{\{u>-\infty\}\cap\Om}[-\chi(u)]H_m(u),$$ it implies that
	\begin{equation}\label{eq3.12}
		\ind_{\Om}	[-\chi(\u)]H_m(\u)\geq \ind_{\{u>-\infty\}\cap\Om}[-\chi(u)]H_m(u).
	\end{equation}
	On the other hand, we also have
	\begin{align*}
		-\chi(\gamma)H_m(\gamma)&=\ind_{\Om}[-\chi(\gamma)]H_m(\gamma)+\ind_{\O\smallsetminus\Om}[-\chi(\gamma)]H_m(\gamma)\\
		&\geq \ind_{\Om}[-\chi(\gamma)]H_m(\gamma).
	\end{align*}
	Note that $\gamma=\u$ on $\Om.$ Therefore, we have
	$$ \ind_{\Om}[-\chi(\gamma)]H_m(\gamma)=\ind_{\Om}[-\chi(\u)]H_m(\u).$$
	This implies that
	\begin{equation}\label{eq3.13}
		-\chi(\gamma)H_m(\gamma)\geq \ind_{\Om}[-\chi(\u)]H_m(\u).
	\end{equation}
	Combining inequality $\eqref{eq3.12}$ and inequality $\eqref{eq3.13}$ we get inequality $\eqref{eq3.15}.$
	\n Thus, it follows from inequality $\eqref{eq3.12a}$ and inequality $\eqref{eq3.15}$ and the construction of $\u,$ we obtain that $\gamma\leq\u$ on $\O.$ This implies that $\psi\leq\u$ on $K.$ Thus, we get the claim $\u\in MSH_m(\O\smallsetminus\overline{\Om})$ as desired.\\
	The claim implies that $H_m(\u)=0$ on $\O\smallsetminus\overline{\Om}.$ Therefore, we obtain that
	\begin{equation}\label{eq3.18}
		\int_{\O}H_m(\u)=\int_{\overline{\Omega}}H_m(\u)
	\end{equation}
	On the other hand, according to equality $\eqref{eq3.12b}$ and the assumption $\eqref{eq3.8a}$, we have 
	$$\int_{\O}[-\chi(\u)]H_m(\u)=\int_{\Om}[-\chi(u)]H_m(u)<+\infty.$$
	Since supp$H_m(\u)\subset\overline{\Om}$, we have 
	\begin{equation}\label{eq3.19}\int_{\overline{\Om}}[-\chi(\u)]H_m(\u)<+\infty.
	\end{equation}
	 Since $\u=0$ on $\partial\O$ and $\Omega\Subset\O,$ we obtain $-\infty\leq\u\leq C<0$  on $\overline{\Om}.$ Note that $\chi$ is a nondecreasing function and $\chi(t)<0$ for all $t<0$, we get $$-1=\chi(-\infty)\leq \chi(\u)\leq\chi(C)<0.$$  This implies that $1\geq-\chi(\u)\geq-\chi(C)>0.$ Combining this with inequality $\eqref{eq3.19}$ we deduce that
	\begin{equation}\label{eq3.20}
		\int_{\overline{\Om} }H_m(\u)<+\infty.
	\end{equation}
	It follows from equality $\eqref{eq3.18}$ and inequality $\eqref{eq3.20}$ that $\int_{\O}H_m(\u)<+\infty.$ Theorem 4.8 in \cite{T19} implies that $\u\in\mathcal{F}_m(\tilde{\u},\O)$, where the function $\tilde{\u}$ defined by $\tilde{\u}=(\lim\limits_{j\to\infty}\u^j)^*$ as in the definition of the class $\mathcal{N}_m$ in the Subsection 2.6 of \cite{Pjmaa}. Specifically, there exists a function $\theta\in\F_m(\O)$ such that 
	\begin{equation}\label{e3.24}\tilde{\u}+\theta\leq\u\leq \tilde{\u}
	\end{equation}
	and 
	$$\int_{\O}H_m(\theta)\leq\int_{\O}H_m(\u)<+\infty.$$
	Since $0\leq -\chi(\theta)\leq 1,$ so we have
	$$\int_{\O}[-\chi(\theta)]H_m(\theta)\leq\int_{\O}H_m(\theta)<+\infty.$$
	According to Proposition 3.14 in \cite{DT23} we have $\theta\in\mathcal{E}_{m,\chi}(\O).$
	Note that since $\u\in\F_m(g,\O),$ there exist $\eta\in\F_m(\O)$ such that $\u\geq g+\eta.$ This implies that $\tilde{\u}\geq \tilde{g}+\tilde{\eta}=\tilde{g}\geq g.$ It follows from inequality $\eqref{e3.24}$  that $\u\geq g+\theta.$ This means we have $\u\in\mathcal{E}_{m,\chi}(g,\O).$ The proof of Case 1 is complete. \\
	{\bf Case 2.} When $\chi(-\infty)=-\infty.$
	By Theorem 3.7 in \cite{DT23}, we infer that $\varphi\in\N^a_m(\Om).$ Note that, we also have $f\in MSH_m(\Om)$ and $u\geq f+\varphi.$ For every $m$-polar sets $A\subset\Om,$ Proposition \ref{pro5.2HP17} implies that
	$$\int_{A}H_m(u)\leq \int_{A}H_m(f+\varphi)\leq 0,$$
	where the last inequality is due to Lemma 5.6 in \cite{HP17}.
	This means that $\ind_{\Om}[-\chi(u)]H_m(u)$ is a positive Radon measure which vanishes on all $m$-polar sets on $\O.$ Moreover, we have 
	$$\int_{\O}\ind_{\Om}[-\chi(u)]H_m(u)=\int_{\Om}[-\chi(u)]H_m(u)<+\infty.$$ By Theorem  3.1 in \cite{PD2024}, there exists  a functions $\u\in\E_{m,\chi}(g,\O)$ such that 
	$$-\chi(\u)H_m(\u)=\ind_{\Om}[-\chi(u)H_m(u)]\,\,\text{on}\,\, \O.$$
	So it remain to prove that $\u\leq u$ on $\Om.$ Indeed, put $w=\max(\u,u)$ on $\Om.$ By Lemma 3.1 in \cite{HQ21}, we have $-\chi(w)H_m(w)\geq-\chi(u)H_m(u)$ on $\Om.$ Therefore, by Theorem 3.8 in \cite{PDjmaa} we infer that $w\leq u.$ It follows that $\u\leq u$ on $\Om.$ So we get desired. The proof is complete.\\
\end{proof}

\section{Approximation of $m$-subharmonic functions in weighted energy classes with given boundary values}
\n Let $\chi: \mathbb{R}^-\to\mathbb{R}^-$ be a nondecreasing continuous function such that $\chi(t)<0$ for all $t<0$ and $\chi\in C^1(\R)$. In this last Section, we will use the result of Section 3 to study approximation of $m$-subharmonic functions in weighted energy classes with given boundary values. 
\begin{theorem}\label{mainth}
	Let $\Omega\Subset \Omega_{j+1}\Subset\Om_j$ are  bounded $m$-hyperconvex domains such that for every compact subset $K$ of $\Om$ we have 
	\begin{equation}\label{e5.1}\lim\limits_{j\to\infty} Cap_m(K,\Om_j)= Cap_m(K,\Om).
	\end{equation}  Assume that $g\in \E_m(\Om_1)\cap MSH_m(\Om_1).$ Then for every function $u\in\mathcal{E}_{m,\chi}(g|_{\Om},\Om)$ satisfying $\int_{\Om}[-\chi(u)+1]H_m(u)<\infty,$ there exists an increasing sequence of functions $u_j\in\mathcal{E}_{m,\chi}(g|_{\Om_j},\Om_j)$ such that $\lim\limits_{j\to\infty}u_j=u$ a.e. on $\Om.$
\end{theorem}
Note that, in \cite{A14}, the author assumes $g\in C(\overline{\Omega})$ while in our result, the condition $g\in C(\overline{\Omega})$ is not necessary.\\
\n We need the following lemma which was essentially proved in \cite{PDjmaa}.
\begin{lemma}\label{lm4.2}
 Let $u\in\mathcal{N}_{m}(f), v\in\mathcal{E}_{m}(\Omega),\,v\leq f$ be such that $-\chi(u)H_{m}(u)\leq -\chi(v)H_{m}(v).$ 
	Assume also that $H_m (u)$ puts no mass on all $m$-polar sets in $\Omega.$
	Then we have $u\geq v$ on $\Om.$
	\end{lemma}
\begin{proof}
Repeating the same argument as in the  Theorem 3.8 in \cite{PDjmaa}. Note that the result of this Theorem will not change if we replace the condition $v\in\mathcal{E}_m(f)$  by the conditions $v\in\mathcal{E}_m(\Om)$ and $v\leq f$.
\end{proof}
\begin{proof}[ Proof of Theorem 4.1]

We outline some main points of our proof. The proof is divided into two cases $\chi(-\infty)>-\infty$ and $\chi(-\infty)=-\infty.$ In both cases, using subextension in Section 3, we construct an incresing sequence in the class $\mathcal{E}_{m,\chi}(f_j,\Om_j).$ Moreover, according to the result in \cite{Pjmaa} (see also \cite{AG23}), there exists an increasing sequence $w_j\in\mathcal{F}_m(f_j,\Omega_j)$ such that $\lim\limits_{j\to\infty}w_j=u$ a.e. on $\Omega.$ Taking the maximum of the two sequences above, we obtain the desired sequence.
Now, we give a detail proof.\\
	We put $g|_{\Om}=f$ and $g|_{\Om_j}=f_j.$ 
	We split the proof into two cases\\
	
	\n	{\bf Case 1.} $\chi(-\infty)>-\infty.$ Without loss of generality, we can assume that $\chi(-\infty)=-1.$\\
	According to Step 1 of Case 1 in the proof of Theorem \ref{th3.3}, there exist $v\in\F_m(f,\Om)$ such that $v\leq u$ and $H_m(v)=\ind_{\{u=-\infty\}}H_m(u).$ Put
	$$v_j=(\sup\{\xi\in\mathcal{E}_m(f_j,\Om_j):\xi\leq v\,\,\text{on}\,\,\Om\})^*.$$
	Repeating the argument as in Step 2 of Case 1 in the proof of Theorem \ref{th3.3}, we obtain that $v_j\in\mathcal{F}_m(f_j,\Omega_j), v_j\leq v$ on $\Omega$ and $H_m(v_j)=\ind_{\Omega}H_m(v).$\\
	We consider the function $h_j$ defined by
	$$h_j=\big(\sup\{\xi\in \mathcal{E}_{m}(\Omega_j): \xi\leq v_j \,\,\text{on}\,\, \Omega_j \,\,\text{and}\,\, -\chi(\xi)H_m(\xi)\geq\alpha \}\big)^*,$$ where
	$$\alpha =\ind_{\{u>-\infty\}\cap\Om}[-\chi(u)]H_m(u)$$ 
	According to the proof of Theorem $\ref{th3.3}$ ( Step 3, Step 4 and Step 5 of Case 1), we have $h_j\in\mathcal{E}_{m,\chi}(f_j,\Om_j).$ 
	From the definition of $h_j$, we see that $[h_j]$ is an increasing sequence. 
	\n Obviously, since $h_j\leq u$ on $\Om,$ we infer that 
	\begin{equation*}
		(\lim\limits_{j\to\infty}h_j)^*\leq u\,\,\text{ on}\,\, \Om.
	\end{equation*}
	On the other hand, we also have $u\in\mathcal{E}_{m,\chi}(f,\Om)\subset\N_m(f,\Om).$ Moreover, it follows from the assumption $\int_{\Om}H_m(u)<+\infty$ and Lemma $\ref{l3.1}$ that $u\in\F_m(f,\Om).$ By the assumption $\ref{e5.1},$ we may apply Theorem 4.5 in \cite{Pjmaa} (see also Theorem 4.1 in \cite{AG23}) to get an increasing sequence of functions $w_j\in\F_m(f_j,\Om_j)$ such that $\lim\limits_{j\to\infty}w_j=u$ a.e. on $\Om.$\\
	Now, we put $u_j=\max(h_j,w_j)$ then we have $u_j \in\mathcal{E}_{m,\chi}(f_j,\Om_j)$ and $[u_j]$ is an increasing sequence such that $\lim\limits_{j\to\infty}w_j=u$ a.e on $\Om.$ \\
	\n	{\bf Case 2.} $\chi(-\infty)=-\infty.$ \n For each $j\in\mathbb{N}^*,$  according to the Case 2 in the proof of Theorem \ref{th3.3}, there exists $v_j\in\mathcal{E}_{m,\chi}(f_j,\Om_j), v_j\leq u\,\,\text{on}\,\,\Om$  and
	\begin{equation*}\label{e4.5}
		-\chi(v_j)	H_m(v_j)= \ind_{\Om}[-\chi(u)]H_m(u)\,\,\text{on}\,\,\Om_j.
	\end{equation*} 
We will apply Lemma 4.2 to prove that $v_j\leq v_{j+1}$ on $\Om_{j+1}.$ Indeed, we will verify the conditions to apply Lemma 4.2.\\
\n (i)	On $\Om_{j+1}$ we have
	$$-\chi(v_j)	H_m(v_j)= \ind_{\Om}[-\chi(u)]H_m(u)= -\chi(v_{j+1})H_m(v_{j+1}).$$
\n (ii)	Obviously, we have $v_j\in \mathcal{E}_{m,\chi}(f_j,\Om_j)\subset \mathcal{E}_{m}(\Om_{j})\subset\mathcal{E}_{m}(\Om_{j+1}).$ We also have $v_{j+1}\in\mathcal{E}_{m,\chi}(f_{j+1},\Om_{j+1})\subset\mathcal{N}_{m}(f_{j+1},\Om_{j+1}).$\\
\n (iii)	Moreover, since $v_{j+1}\in\mathcal{E}_{m,\chi}(f_{j+1},\Om_{j+1}),$ there exists a function $\varphi\in\mathcal{E}_{m,\chi}(\Om_{j+1})$ such that $v_{j+1}\geq f_{j+1}+\varphi.$ Thus, for every $m$-polar sets $A\subset\Om_{j+1},$ Proposition \ref{pro5.2HP17} implies that
	\begin{equation}\label{e4.2}\int_{A}H_m(v_{j+1})\leq \int_{A}H_m(f_{j+1}+\varphi).
	\end{equation}
	Since $\chi(-\infty)=-\infty,$ by Theorem 3.7 in \cite{DT23} we have $\varphi\in\mathcal{E}_{m,\chi}(\Om_{j+1})\subset \mathcal{E}_m^a(\Om_{j+1}).$ Note that we have $f_{j+1}\in MSH_m(\Om_{j+1}).$ Therefore, by Lemma 5.6 in \cite{HP17} we get
	\begin{equation}\label{e4.3}
		\int_{A}H_m(f_{j+1}+\varphi)=0.
	\end{equation}  
	Combining inequality \eqref{e4.2} and equality \eqref{e4.3} we have $H_m(v_{j+1})$ vanishes on all $m$-polar sets on $\Om_{j+1}.$ \\
	It follows from (i), (ii) and (iii) that on $\Omega_{j+1},$ we have $-\chi(v_j)	H_m(v_j)=-\chi(v_{j+1})H_m(v_{j+1})$ with $v_j\in\mathcal{E}_{m}(\Om_{j+1}), v_j\leq f_{j+1},$ $v_{j+1}\in\mathcal{N}_{m}(f_{j+1},\Om_{j+1})$  and $H_m(v_{j+1})$ vanishes on all $m$-polar sets on $\Om_{j+1}.$ According to Lemma \ref{lm4.2}, we obtain that $v_j\leq v_{j+1}$ on $\Om_{j+1}.$ This mean that we have $[v_j]$ is an increasing sequence. Since $v_j\leq u$ on $\Om,$ we have $(\lim\limits_{j\to\infty}v_j)^*\leq u.$\\
	Repeating the argument as in the last part in the proof  of the Case 1, we get desired result.
	The proof is complete.
\end{proof}
\section*{Declarations}
\subsection*{Ethical Approval}
This declaration is not applicable.
\subsection*{Competing interests}
The authors have no conflicts of interest to declare that are relevant to the content of this article.
\subsection*{Availability of data and materials}
This declaration is not applicable.

\end{document}